\documentclass[11pt]{article}
\newcommand{\version}{July 2008}

\usepackage{amsthm,amsfonts,amsmath}
\usepackage{latexsym}
\swapnumbers
\theoremstyle{plain}

\newtheorem{theorem}{THEOREM}[section]
\newtheorem{lm}[theorem]{LEMMA}
\newtheorem{cl}[theorem]{COROLLARY}

\theoremstyle{definition}

\theoremstyle{remark}
\newcommand{\upchi}{\raise1pt\hbox{$\chi$}}
\newcommand{\R}{{\mathord{\mathbb R}}}

\newcommand{\F}{{\mathcal{F}}}

\renewcommand{\|}{{\Vert}}
\numberwithin{equation}{section}
\def\MM{{\cal M}}
\def\F{{\cal F}}
\def\d{{\rm d}}
\def\G{{\cal G}}

\begin{document}

\markboth{\scriptsize{CCELM \version}}{\scriptsize{CCELM \version}}

\title{\bf{Displacement convexity and minimal fronts at phase boundaries}}
\author{\vspace{5pt} E. A.  Carlen$^1$, M. C. Carvalho$^2$, R. Esposito$^3$, J.L. Lebowitz$^4$
 and R. Marra$^5$\\}

\date{\version}
\maketitle

\def\O{\Omega}

\footnotetext                                                                         
[1]{ Department of Mathematics, Rutgers University,
Piscataway, NJ  08854, U.S.A. Work partially
supported by U.S. National Science Foundation
grant DMS 06-00037.  }
\footnotetext 
[2]{Department of Mathematics and CMAF, University of Lisbon,
1649-003 Lisbon, Portugal. Work partially supported by POCI/MAT/61931/2004}
\footnotetext 
[3]{Dip. di
Matematica, Universit\`a di L'Aquila, Coppito, 67100 AQ, Italy}
\footnotetext 
[4]{Departments of Mathematics and
Physics, Rutgers University, New Brunswick, NJ
08903, U.S.A.}
\footnotetext 
[5]{Dipartimento di Fisica and Unit\`a INFN, Universit\`a di
Roma Tor Vergata, 00133 Roma, Italy. \\
\copyright\, 2007 by the authors. This paper may be reproduced, in its
entirety, for non-commercial purposes.}
                                                                      
\begin{abstract}

We show that certain free energy functionals that are not convex with respect to the usual
convex structure on their domain of definition, are strictly convex in the sense of displacement convexity under a natural change of variables. 
We use this to show that in certain cases, the only critical points of these functionals are
minimizers. This approach based on displacement convexity permits us to treat
multicomponent systems as well as single component systems. The developments produce
new examples of displacement convex functionals, and, in the multi-component setting, jointly displacement convex functionals.

\end{abstract}

\bigskip
\centerline{Mathematics Subject Classification Numbers:  49S05, 52A40, 82B26}

\section{Introduction} \label{intro}

\subsection{The variational problem} \label{intro1}

\medskip
We consider  minimization problems for a type of functional that arises in the study of phase segregation in statistical mechanical systems.  Let $F(m)$ be a  function on the real line
that is continuous and strictly positive except at $m= a$ and $m=b$ with $a< b$. A good example to bear in mind is the ``double well potential'' 
$$F(m) = \frac{1}{4}(m^2 -1)^2\ ,$$
where of course $a=-1$ and $b=1$.

Let ${\cal C}_{a,b}$ be the set of  measurable functions $m(x)$ from $\R$ to $\R$ such that
(for some representative)
$$\lim_{x\to -\infty} m(x) =a\qquad{\rm and}\qquad  \lim_{x\to +\infty} m(x) =b\ .$$
The numbers $a$ and $b$ represent the values of the order parameter $m$ in two phases of a statistical mechanical system. For example, $m=a$ might correspond to a vapor phase, and $m=b$ to a liquid phase. 

A
function $m(x)$ in  ${\cal C}_{a,b}$ denotes a  possible one-dimensional {\it transition profile} across the boundary segregating  the two different phases. The actual profile that one would expect to see would be one that minimizes the free energy cost of making such a transition. 
The free energy functional  $\F$ to be minimized on  ${\cal C}_{a,b}$ will in some cases of interest have the form, c.f. \cite{dopt1}, 
\begin{equation}\label{free1}
\F(m) = \int_\R F(m(x))\d x + \frac{1}{2}\int_\R\int_\R (m(x) - m(y))^2J(x-y)\d x\d y\ ,
\end{equation}
where $J(x)$ is a non-negative integrable  function on $\R$. 

The term $ \int_\R F(m(x))\d x$
is due to short range interactions and entropy effects and is normalized so that it vanishes in the pure phases, when $m(x)=a$ or $m(x)=b$, 
while the term $\int_\R\int_\R (m(x) - m(y))^2J(x-y)\d x\d y$ is due to
long range interactions. This long range term in the free energy suppresses sharp transitions,
 as does the gradient term in the familiar but purely phenomenological  Van der Waals model \cite{VW}.  For more discussion of the physical context of the problem, see \cite{CCELM2}.

Much useful information can be deduced from the specific form of the minimizing profiles. 
In particular, the surface tension at a two dimensional phase boundary in physical three dimensional space is the minimum value of $\F(m)$ on
${\cal C}_{a,b}$; see \cite{alb2} and Section~\ref{binsec}  for more information. Hence we ask:

\medskip

\noindent{$\bullet$} {\it What is the minimum value of $\F(m)$ as $m$ ranges over 
${\cal C}_{a,b}$, and  are the minimizing profiles, if any, unique up to translation?}

\medskip

Actually, the existence of minimizers is relatively simple to prove using the rearrangement inequalities to be discussed below.  However, because of the translation invariance, they are never unique: Any translate of a minimizer is again a minimizer. It is less  simple to show that this is the only degeneracy. 

\subsection{Displacement convexity  and uniqueness of fronts} \label{intro2}

For a particular choice of $F$ in the free energy functional specified in (\ref{free1}), the minimizing profile  problem  has been solved in a series of papers \cite{dopt1},\cite{dopt2} by De Masi, Orlandi, Triolo and Presutti, building on previous unpublished work of Dal Passo and 
de Mottoni \cite{dpdm}
Their solution involves the construction of a dynamics that is dissipative for the free energy functional, and then a careful analysis of limits along the
time evolution for this dynamics.

Another approach that we further develop here has been introduced by Alberti and Bellettini
\cite{ab}, \cite{alb}.  They discovered an alternative convex structure which renders the
variational problem for (\ref{free1}) convex, and used this to study the existence problem in 
\cite{ab}. Later, Alberti  \cite{alb} returned to the problem and proved a uniqueness result that affirmatively answers the question raised above for this one component model.

Our goal here is to treat certain two component systems.  Motivated by this problem, we were led
to reconsider the single component problem from the point of view of 
McCann's notion of {\it displacement convexity} \cite{McC1}.  In fact, the minimization problem for  (\ref{free1}) is 
challenging largely because the functional ${\cal F}$ is 
{\it not convex} on ${\cal C}_{a,b}$ in the usual way:
For $0 < \lambda  < 1$, and $m_0$ and $m_1$
in  ${\cal C}_{a,b}$, define $m_\lambda = (1-\lambda)m_0 + \lambda m_1$ and note that $m_\lambda \in  {\cal C}_{a,b}$. However,  due
to the non convexity of the potential function $F$, it is {\it not}  true in general that  ${\cal F}(m_\lambda) \le (1-\lambda){\cal F}(m_0) +  \lambda{\cal F}(m_0)$.

In \cite{McC1}, McCann, building on groundbreaking work of Brenier \cite{bren}, introduced an alternative convex  structure on the space of probability densities on $\R^n$, and used this
to prove existence and uniqueness results for minimizers of functionals that were not convex in the usual sense. We shall show here that the minimization problem for (\ref{free1}), as well as for
a two component model of this type, can be handled within this framework. In the process, we
provide two new examples of strictly displacement convex functionals, the second of which is jointly displacement convex.   It turns out that the alternative convex structure introduced in \cite{ab}
is equivalent to the displacement convexity in this one dimensional setting, although the approach is quite different. We shall see that developing the alternative convex structure explicitly in terms of displacement convexity has advantages, especially for the two component system, when one
seeks to prove a uniqueness result. Moreover, as we show in the final section, our results for the two component system may be applied to the single component system in higher dimensions, yielding a new uniqueness theorem for monotone solutions of the Euler-Lagrange equation.

We now describe the alternate convex structure with respect to which $\F$ {\it is} convex.  This second convex structure cannot be defined on all of  ${\cal C}_{a,b}$, but only on the subset ${\cal M}_{a,b}$ 
consisting of {\it right continuous monotone}
profiles. Nothing is lost in this restriction, as  rearrangement inequalities 
show that minimizers of $\F$ on   ${\cal C}_{a,b}$ must actually be monotone, so that they have a right continuous  version belonging  ${\cal M}_{a,b}$;
see \cite{alb} and Theorem \ref{rear} below.

Any right--continuous profile $m(x)$  in ${\cal M}_{a,b}$
can be written in the form
\begin{equation}m(x) = a + (b-a)\int_{(-\infty,x]} {\rm d}\mu(y)\label{eq12}\end{equation}
where $\mu$ is a uniquely determined probability measure on $\R$. This identification  of  ${\cal M}_{a,b}$
and the set of probability measures on $\R$
allows us to look at $\F$ as a functional defined on probability measures. 

This is a useful perspective since there is an  alternative convex structure on the set of probability measures on $\R$
(or more general domains) that was introduced by McCann, and which we describe below.
A functional on probability measures
is said to be {\em displacement convex} if it is convex with respect to this alternative structure. We shall show here that $\F$, regarded as a functional  on probability measures is, in fact, displacement convex.
Using this, we shall  show that any solution in ${\cal M}_{a,b}$  of the Euler--Lagrange  equation for the variational problem concerning  (\ref{free1}) 
\begin{equation}\label{el}
m(x) = \frac{1}{\widehat J}\left(\int_\R J(x-y)m(y)\d y - F'(m(x))\right)\ ,
\end{equation}
where 
\begin{equation}\label{toma}
\widehat J = \int_\R J(x)\d x\ , 
\end{equation}
is in fact a minimizer. 
Solutions to (\ref{el}) can easily be constructed by iteration and using these surface tensions may be readily computed.

This solution to the variational problem has the advantage of applying also to free energy functionals in certain 
multicomponent systems, in which the determination of the minimizers has not been previously treated. Indeed,  our motivation was to be able to rigorously determine the surface tension in such systems.   However,  we shall first present 
our simple solution of the minimization problem for the single component free energy
functional $\F$  specified in (\ref{free1}), and then treat the multicomponent case.

\section{The alternative convex structure}\label{displace}

\subsection{The reduction to monotone profiles} %
\medskip

First of all, notice that if we seek to minimize $\F$ on ${\cal C}_{a,b}$, we need only consider profiles
$m$ for which $a\le m(x) \le b$ for all $x$.  Indeed, for any $m\in {\cal C}_{a,b}$, define $\widehat m$
by
$$\widehat m(x) = \min\{b , \max\{ a,m(x)\}\}\ .$$
Then  $\F(\widehat m) \le \F(m)$ with equality only in case $\widehat m = m$, since otherwise replacing $m$ by $\widehat m$ lowers both the potential and the interaction terms.

We now recall a notion of rearrangement due to Alberti \cite{alb}. For any Borel measurable set $A$, 
let $|A|$ denote its Lebesgue measure. The rearrangement is defined  for  Borel sets $A\subset \R$ such that $|A\Delta(0,\infty)| < \infty$, where $A\Delta B = A\backslash B \cup B\backslash A$ is the symmetric
difference of $A$ and $B$.
For such a set $A$, define the rearranged set $A^*$ by
$$A^* = [\alpha,\infty)\qquad{\rm where}\qquad \alpha =|(0,\infty)\backslash A| -  |A\backslash (0,\infty)|\ .$$

Any function $m$ in ${\cal C}_{a,b}$ that takes values in $[a,b]$ can be represented in ``layer--cake'' form:
$$m(x) = \int_a^b 1_{\{m > z\}}(x){\rm d}z+ a\ .$$
For each $z\in (a,b)$, the set  $\{m > z\}$ certainly has the property that $|\{m > z\}\Delta (0,\infty)| <\infty$.
Hence one can define the rearrangement of $m$ itself through
$m^*(x) =  \int_a^b \left(1_{\{m > z\}}\right)^*(x){\rm d}z+a$. (Applying the rearrangement to a monotone increasing function, one simply obtains the right continuous version.)

Alberti shows that for any two such functions $m_1$ and $m_2$, 
$$\int_{\R}|m_1^*(x)- m_2^*(x)|^2{\rm d}x \le  \int_{\R}|m_1(x)- m_2(x)|^2{\rm d}x\ .$$
In particular, with $m$ being any function in ${\cal C}_{a,b}$ that takes values in $[a,b]$, and $h$ any real number, let
 $m_1(x) = m(x)$, and $m_2(x) = m(x+h)$. Then
 $$\int_{\R}|m^*(x)- m^*(x+h)|^2{\rm d}x \le  \int_{\R}|m(x)- m(x+h)|^2{\rm d}x\ ,$$
 so that
 $$\int_\R\left( \int_{\R}|m^*(x)- m^*(x+h)|^2{\rm d}x\right)J(h){\rm d}h \le$$ $$  
 \int_\R\left(\int_{\R}|m(x)- m(x+h)|^2{\rm d}x  \right)J(h){\rm d}h\ .$$
 This of course means that
 \begin{equation}\label{al1}
  \int_\R\int_\R (m^*(x) - m^*(y))^2J(x-y)\d x\d y \le  \int_\R\int_\R (m(x) - m(y))^2J(x-y)\d x\d y\ .
  \end{equation}
 In fact, Alberti shows (see Theorem 2.11 in \cite{alb}) that there is equality in (\ref{al1}) if and only if $m = m^*$.

Of course, $\int_\R F(m^*(x))\d x = \int_\R   F(m(x))\d x$, and so we have
$\F(m^*) \le \F(m)$ with equality if and only if $m = m^*$. Thus, we may restrict our search for minimizers to ${\cal C}_{a,b}$, the subset of monotone increasing profiles in ${\cal C}_{a,b}$.

\subsection{Displacement convexity of $m\mapsto \int_\R  F(m(x))\d x$}
\medskip

As noted in (\ref{eq12}), if $m$ is any profile in 
${\cal M}_{a,b}$, then $(m(x) -a)/(b-a)$ is the cumulative distribution function of a uniquely determined
probability measure $\mu$:
$$\frac{m(x) -a}{b-a} = \int_{(-\infty,x]}  \d \mu(y)\ .$$

For each $m$ in $\MM_{a,b}$, define $x(m)$ to be the inverse function: For $m\in (a,b)$,
\begin{equation}\label{6}
x(m) = \inf\{\ x\ :\ m(x) > m\ \}\ .
\end{equation}
Then of course, $m(x)$ is the inverse function of $x(m)$, so that for $x$ in $\R$,
\begin{equation}\label{7}
m(x)  = \inf\{\ m\ :\ x(m) > x\ \}\ .
\end{equation}
Let ${\rm d}x(m)$ denote  the Lebesgue-Stieltjes measure on $[a,b]$ induced by the monotone
function $x(m)$. (In the terminology introduced below,  ${\rm d}x(m)$ is
 the push-forward of Lebesgue measure on $\R$ under $m$.)   Then one can rewrite
$$\int_\R F(m(x))\d x = \int_a^bF(m)\d x(m)\ .$$
Let  $m_0$ and $m_1$ be any two elements of $\MM_{a,b}$,
and let $x_0$ and $x_1$ denote their respective inverse functions. Then for any  $\lambda \in (0,1)$,
define $x_\lambda(m)$ by
\begin{equation}\label{8}
x_\lambda(m) = (1-\lambda) x_0(m) + \lambda x_1(m)\ .
\end{equation}
Note that  $x_\lambda$ is also the inverse function of an element of $\MM_{a,b}$, which we shall call $m_\lambda$.
That is,
\begin{equation}\label{goodint}
m_\lambda(x) = \inf\{\ m\ :\  (1-\lambda) x_0(m) + \lambda x_1(m) > x\ \}\ .
\end{equation}
Note that ${\rm d}x_\lambda(m)$, the  the Lebesgue-Stieltjes measure on $[0,1]$ induced by the monotone
function $x_\lambda(m)$, satisfies ${\rm d}x_\lambda(m) = (1-\lambda){\rm d}x_0 + \lambda {\rm d}x_1$.
Then,
\begin{eqnarray}\label{9}
\int_\R F(m_\lambda(x))\d x &=& \int_a^b F(m)\d x_\lambda(m)\nonumber\\
&=&(1- \lambda) \int_a^b F(m)\d x_0(m) + \lambda
 \int_a^b F(m)\d x_1(m)\nonumber\\
&=& (1-\lambda) \int_\R F(m_0(x))\d x + \lambda  \int_\R F(m_1(x))\d x\ .\nonumber\\
\end{eqnarray}

This tells us that along the interpolation $m_\lambda$ between $m_0$ and $m_1$ provided by 
(\ref{goodint}), the function $\lambda \mapsto  \int_\R F(m_\lambda(x))\d x$ is affine, and in particular, is convex.
This is not the case for the standard interpolation given by
\begin{equation}\label{usint}
\widetilde m_\lambda(x) = (1-\lambda)m_0(x) + \lambda m_1(x)\ ,
\end{equation}
since 
$\lambda \mapsto  \int_\R F(\widetilde m_\lambda(x))\d x$
is {\it not}, in general, convex.  That is,  taking convex combinations in terms of the inverse function $x(m)$, as in (\ref{goodint}), instead of $m(x)$ itself, as in (\ref{usint}), 
has ``cured'' the non-convexity of the  functional $m \mapsto \int_\R F(m(x))\d x $.

Of course, this will only be useful if the functional
\begin{equation}\label{5}
m \mapsto \int_\R\int_\R (m(x) - m(y))^2J(x-y)\d x\d y\ ,
\end{equation}
which was convex in the usual way, is still convex with the new convex structure. This is not at all obvious, but the main result of the next section asserts that this is the case. 

The approach of Alberti and Bellettini \cite{ab}, which we discovered only after our work was complete, was
to rewrite the interaction directly in terms of $x_m$, and to show that it is convex. 

However, it turns out that the convex structure in (\ref{8})  is something that is by now
well--known; it is the {\it displacement convexity}
structure introduced by McCann. Making this connection will facilitate showing the {\it strict} convexity of
$m \mapsto \int_\R\int_\R (m(x) - m(y))^2J(x-y)\d x\d y$ under this convex structure.
This point was left open in \cite{ab}, who explicitly asked whether one could extend the
ideas to give a direct proof of uniqueness. Although Alberti \cite{alb} did later return to address the issue, we shall see here that the strict convexity is quite clear from the perspective of displacement convexity.

Displacement convexity is usually introduced as a convex structure in a set of probability measures. 
Given a probability measure $\mu_0$ on $\R$, and a measurable map $T:\R\to \R$, we define the {\it push forward of $\mu_0$ under $T$},
$T\#\mu_0$, by
\begin{equation}\label{10}
\int_\R \phi(T(x))\d \mu_0(x) = \int_\R \phi(y)\d(T\#\mu_0)(y)\ ,
\end{equation}
for all bounded, continuous  functions $\phi$.

Given two probability measures $\mu_0$ and $\mu_1$ on $\R$, there  is a unique {\it monotone} map $T$
such that $T\#\mu_0 = \mu_1$. To see what it must be, fix any $a\in R$, let $\phi_a$
be the step function 
$\phi_a(x) = 1_{(-\infty,a]}(x)$.
Then, by definition, we must have
${\displaystyle \int_\R \phi_a(T(x))\d \mu_0(x) = \int_\R \phi_a(y)\d \mu_1(y)}$,
and hence
\begin{equation}\label{11}
\int_{-\infty}^{T^{-1}(a)}\d \mu_0 =  \int_{-\infty}^{a}\d \mu_1\ .
\end{equation}
Let $m_0$ and $m_1$ be the cumulative distribution functions of $\mu_0$ and $\mu_1$, respectively.
Then (\ref{11}) entails that $m_0(T^{-1}(a)) = m_1(a)$ for all $a$, or, what is the same thing
\begin{equation}\label{12}
m_0(a) = m_1(T(a))\ 
\end{equation}
for all $a$. 
As long as $m_1$ is free of ``flat spots'', so that the inverse function does the expected thing,
this leads to 
\begin{equation}\label{13}
T(a) = x_1(m_0(a))\ .
\end{equation}
As long as $\mu_0$ and $\mu_1$ have strictly positive densities, (\ref{13}) does indeed define a monotone map $T$,
and then it is very easy to see that with $T$ defined by (\ref{13}),  $T\#\mu_0 = \mu_1$, and in fact, this is true 
without further technical hypotheses; see \cite{V} for more information.

We now interpolate the map $T$, and hence the corresponding probability measures 
$\mu_0$ and $\mu_1$ and the corresponding cumulative distribution functions $m_0$ and $m_1$ as well.
For all $\lambda \in [0,1]$, define $T_\lambda$ by
\begin{equation}\label{14a}
T_\lambda(x)  = (1-\lambda)x + \lambda T(x)\ .
\end{equation}
If we define $x_\lambda(m)$ by
$$x_\lambda(m) = T_\lambda(x_0(m))\ ,$$
then clearly $x_\lambda$ is given by (\ref{8}). 

The displacement convex structure on probability measures on $\R$ is given by 
$\mu_\lambda = T_\lambda\#\mu_0$,
and so it is nothing other than the convex structure (\ref{8}), expressed in terms of probability measures instead of
cumulative distribution functions.  When $\mu_1$ and $\mu_2$ have strictly positive densities, so that $T$ is given by (\ref{12}),
we denote the density of $\mu_\lambda$ by $\rho_\lambda$, and write
\begin{equation}\label{15}\rho_\lambda = T_\lambda\#\rho_1\ .
\end{equation}
We summarize the main result of this section in a theorem:

\medskip
\begin{theorem}\label{affine} Let $\lambda \mapsto m_\lambda$ be the displacement interpolation between $m_0$ and $m_1$
 in  $\MM_{a,b}$. Then for $0 \le \lambda \le 1$,
 $$\int_\R F(m_\lambda(x)){\rm d}x = (1-\lambda) \int_\R F(m_0(x)){\rm d}x + \lambda
 \int_\R F(m_1(x)){\rm d}x\ .$$
 \end{theorem}

\section{Displacement convexity of the interaction energy}\label{lower}

Let $\MM$ denote the class of cumulative distribution functions on $\R$.
Making the obvious change of variables, we will assume without loss of generality that $a=0$
and $b=1$ and we will set $\MM_{0,1}=\MM$.

Given any $m\in \MM$, let $\mu$ be the corresponding probability measure, so that
$m(x) = \int_{-\infty}^x\d \mu(y)$.  The first step in the investigation of the interaction energy is to rewrite it as a functional of 
$\mu$ instead of $m$. This is done in the following lemma:

\medskip
\begin{lm}\label{iaconv}   Assume that $\int_\R |s|J(s)\d s < \infty$. Define $W$ in terms of $J$ by setting
\begin{equation}\label{wform}W(u) = \int_{u}^\infty (s - u)(J(s)+ J(-s))\d s\ .\end{equation}
for $u\ge 0$ and $W(u)=W(-u)$ for $u<0$. Then
$$\int_\R\int_\R (m(x) - m(y))^2J(x-y)\d x\d y = \int_\R\int_\R W(z-w)\d \mu(z)\d \mu(w)\ .$$
$W$ is a symmetric function, and  is convex on $(0,\infty)$ and on $(-\infty,0)$, though not on all of $\R$. 
\end{lm}

\medskip
\noindent{\bf Proof:} Since for $x < y$, $m(x) - m(y) = \int_x^y \d \mu(z)$, we have from the Fubini Theorem that
\begin{eqnarray}
&\phantom{*}&\int_\R\int_\R (m(x) - m(y))^2J(x-y)\d x\d y =\nonumber\\
&\phantom{*}&\int_\R\int_\R\left [\int_\R \int_\R1_{[x,y]}(z) 1_{[x,y]}(w)\d \mu(z)\d \mu(w)\right]J(x-y)\d x\d y =\nonumber\\
&\phantom{*}&\int_\R\int_\R\left [\int_\R \int_\R1_{[x,y]}(z) 1_{[x,y]}(w)J(x-y)\d x\d y  \right]\d \mu(z)\d \mu(w)\ .\nonumber\\
\end{eqnarray}
Thus if we define $V(z-w)$ by
${\displaystyle V(z-w) = \int_\R \int_\R1_{[x,y]}(z) 1_{[x,y]}(w)J(x-y)\d x\d y}$
we  have
$$\int_\R\int_\R (m(x) - m(y))^2J(x-y)\d x\d y = \int_\R\int_\R V(z-w)\d \mu(z)\d \mu(w)\ .$$

We next show that $V = W$. 
To do this,  write 
$$J_+(x) = \begin{cases} 
J(x)\, 
&\text{for $x>0$}\  , \\
0 &\text{for $x\leq 0$} \ 
\end{cases} \qquad{\rm and}\qquad J_-(x) = J(x) - J_+(x)\ .$$
We first consider 
$\int_\R \int_\R1_{[x,y]}(z) 1_{[x,y]}(w)J_+(x-y)\d x\d y$.  Make the change of variables $s = y-x$ and $t= (x+y)/2$.
Then $\d x\d y = \d s \d t$, and 
${\displaystyle1_{[x,y]}(z) 1_{[x,y]}(w) = 1_{[t-s/2,t+s/2]}(z) 1_{[t-s/2,t+s/2]}(w)}$.
This quantity is zero unless
$|z-w|\le s$ and $ |2t - (x+w)| \le s - |z-w|$,
in which case it is one. Therefore
\begin{eqnarray}\int_\R \int_\R1_{[x,y]}(z) 1_{[x,y]}(w)J_+(x-y)\d x\d y &=&
\int_{|z-w|}^\infty \left(\int_{(z+w)/2 -(s -|z-w|)/2}^{(z+w)/2 +(s -|z-w|)/2}\d t\right)J_+(s)\d s\\
&=& \int_{|z-w|}^\infty  (s -|z-w|)J_+(s)\d s\ .\nonumber
\end{eqnarray}

Doing the same calculation for the part involving $J_-$, we obtain that $V =W$ where $W$ is given by (\ref{wform}). 
Also, for $u>0$,
$W'(u) = -\int_{u}^\infty(J(s)+ J(-s))\d s$,
and so
$W''(u) = J(u) + J(-u)$,
which is non--negative. Thus, $W$ is convex on $(0,\infty)$, and on $(-\infty,0)$ by symmetry. However, it is not convex on the whole real line. Notice that $W(0) = \int_0^\infty s(J(s)+J(-s))\d s > 0$, while
$\lim_{u\to\pm \infty}W(u) = 0$.
\qed
\medskip

We now prove the main result of this section:
\medskip

\medskip
\begin{theorem}\label{interconv} Let $\lambda \mapsto m_\lambda$ be the displacement interpolation between $m_0$ and $m_1$
 in  $\MM$, as defined in (\ref{goodint}). Then for $0 < \lambda < 1$,
 \begin{eqnarray}\label{conv1}
 \int_\R\int_\R (m_\lambda(x) - m_\lambda(y))^2J(x-y)\d x\d y &\leq&(1-\lambda) \int_\R\int_\R (m_0(x) - m_0(y))^2J(x-y)\d x\d y\nonumber\\
  &+& \lambda
\int_\R\int_\R (m_1(x) - m_1(y))^2J(x-y)\d x\d y\ .\nonumber\\
\end{eqnarray}
If  $J$ is strictly positive on some interval, and $m_0$ has a strictly positive derivative almost everywhere, 
 there is equality if and only if $m_1$ is a translate of $m_0$.
\end{theorem}

 \medskip
 
 \noindent{\bf Proof:} If $W$ were convex on all of $\R$, the displacement convexity of the interaction energy would be
 a classical result of McCann \cite{McC1}.  However, in one dimension, the partial convexity of $W$ that was established in 
 Lemma \ref{iaconv} suffices, as observed by  Blower \cite{b}.  This is because the map 
 $T_\lambda$ is monotone
for all $\lambda$. Therefore, if $z > w$, $T_\lambda(z) > T_\lambda(w)$ for all $\lambda$. Hence, as we vary $\lambda$,
$T_\lambda(z) - T_\lambda(w)$ stays in a domain of convexity of $W$.

Hence from (\ref{8}), if $\d \mu_\lambda = T_\lambda\#\d\mu_0$ is the displacement  interpolation between $\d \mu_0$ and $\d \mu_1$,
\begin{equation}\label{apr1}
\int_\R\int_\R W(z-w)\d \mu_\lambda(z)\d \mu_\lambda(w) = 
\int_\R\int_\R W(T_\lambda(z) - T_\lambda(w))\d \mu_0(z)\d \mu_0(w)\ .
\end{equation}
Define the map 
$S(x)$   by $S(x) = T(x) -x$. Then,
 we can rewrite (\ref{apr1}) as
\begin{equation}\label{apr2}
\int_\R\int_\R W(z-w)\d \mu_\lambda(z)\d \mu_\lambda(w) =
\int_\R\int_\R W([z-w] +\lambda[S(z)-S(w)])\d \mu_0(z)\d \mu_0(w)\ .
\end{equation}

By the remarks made above, the right hand side is clearly a  convex function of $\lambda$.
In fact, under mild assumptions on $\mu_0$ or $J$, it is strictly convex 
unless $T$ is simply a translation.  

To see this {\it formally}, let $J$ be symmetric for simplicity of notation, and
differentiate the right hand side of (\ref{apr2}) twice in $\lambda$, finding
$$ \int_\R\int_\R 2J\left([z-w] + \lambda [S(z) - S(w)]\right)[S(z) - S(w)]^2\d \mu_0(z)\d \mu_0(w)\ .$$
If this vanishes for all $\lambda$, then
${\displaystyle\int_\R\int_\R 2J(z-w)[S(z) - S(w)]^2\d \mu_0(z)\d \mu_0(w)= 0}$.
If  $J$ is strictly positive and  if $\mu_0$ has a strictly positive density, then this is possible if and only if $S$ is constant, and that of course means that $T$ is a translation.

 To make this argument rigorous,  and to relax the hypotheses, let $f(\lambda)$ denote the right hand side of (\ref{conv1}) minus the left hand side.   Then, with $g(z,w,\lambda)$ defined by
 $$g(z,w,\lambda) = [\lambda W(z-w)  + (1-\lambda) W((z-w) + (S(z)-S(w)))] $$ $$- W((z-w) + \lambda(S(z)- S(w)))\ ,$$
 we have
 $f(\lambda) = \int_\R\int_\R  g(z,w,\lambda)\d \mu_0(z)\d \mu_0(w)$.
 Since the integrand is non negative, we have for any measurable subsets $A$ and $B$ of $\R$,
 \begin{equation}\label{conv3}
 f(\lambda) \ge \int_A\int_B  g(z,w,\lambda)\d \mu_0(z)\d \mu_0(w)\ .
 \end{equation}
 
Suppose that $J$ is strictly positive on the open interval $I = (y_0-\delta/2, y_0+\delta/2)$.  Then $I$ is an interval of strict convexity of $W$, so that whenever
 $w-z\in I$, $\lambda \mapsto  g(z,w,\lambda)>0$ on $(0,1)$ unless $S(z) = S(w)$.  However, 
if $S$ is not constant almost everywhere, we can find
an arbitrarily small interval about some $z_0$ on which it has strictly positive oscillation.  In particular, we can find
a $z_0$ and an $\epsilon>0$ so that  $\int_{z_0-\delta/2}^{z_0+\delta/2}(S(z) - c)^2\d z > \epsilon$ for all $c$.
Let $A = (z_0-\delta/2,z_0+\delta/2)$, and let   $B =  (y_0+x_0-\delta/2, y_0+x_0+\delta/2)$.  Then for all $z$ in $A$ and $w$ in $B$, $z-w$ belongs to $I$. Moreover, for every $w$ in $B$, $\int_{A}(S(z) - S(w))^2\d w>0$, so $|S(z) - S(w)|>0$
on a subset of $A$ of positive Lebesgue  measure.  Since $\mu_0$ has a strictly positive density, this ensures that
the right hand side of (\ref{conv3}) is strictly positive.\qed
 
 \medskip
 It is clear that the conditions on $J$ and $m_0$ that are invoked to ensure strict convexity can be relaxed, though they are
 already quite general.

We close this section a remark. If the profile $m$ is continuously differentiable with $m'(x) = \rho(x)$,
and $\int_\R J(x)\d x =1$, 
then
$$\lim_{h\to 0}\int_\R\int_\R \frac{(m(x) - m(y))^2}{h^2}\frac{1}{h}J\left(\frac{x-y}{h}\right)\d x\d y = 
\int_\R \rho^2(x)\d x\ .$$
It is already well known that the functional $\rho \mapsto \int_\R \rho^2(x)\d x$ is displacement convex, so the fact that
Theorem \ref{interconv} gives another proof of this is not of great interest. However, the connection between the two functionals
at least gives one suggestion as to why the interaction functional might be expected to be displacement convex.

  \medskip

\section{For the functional ${\cal F}$, critical points are minimizers}

\medskip

\begin{theorem} If $m_0$ is any critical point of ${\cal F}$ in $\MM$, and  $m$ is any other profile in $\MM$, then
${\cal F}(m) \ge {\cal F}(m_0)$
and there is equality if and only if $m$ is a translate of $m_0$. 
\end{theorem}

\noindent{\bf Proof:} Let $m_\lambda$ be the displacement interpolation between
$m_0$ and $m$. Then $\lambda \mapsto \F(m_\lambda)$ is convex,
and the derivative is zero at $\lambda = 0$. Hence $m_0$ is a minimizer of $\F$, so that
$\F(m_\lambda) \ge \F(m_0)$, and there is equality if and only if $\lambda \mapsto \F(m_\lambda)$
is constant. But in this case, the strict displacement convexity of ${\cal F}$ ensures that $m$ is a translate of $m_0$. 
\qed

\section{Fronts in a binary fluid model}\label{binsec}
\medskip
We now turn to the study of the analogous problem for a binary fluid model. The binary fluid model  has been investigated in \cite{CCELM} and  \cite{CCELM2}, and we refer to those papers for details. Although the arguments apply to that setting in full generality, we  discuss here only a special case where the non local interaction is only between particles of different species and the local term is purely entropic, for the sake of brevity.
For further information and a numerical investigation of the minimizing fronts, see \cite{belm}.

In what follows, $m(x)$ and $n(x)$ represent the particle number densities of two different species of particles contained in some bounded domain  $\Omega$ in $\R^n$. 
Consider the functional
${\cal F}$ defined  by
\begin{eqnarray}
{\cal F}(m,n) &=& \int_\Omega m(x)\ln m(x){\rm d}x  +
\int_\Omega n(x)\ln n(x){\rm d}x \nonumber\\
 &+&\beta  \int_\Omega \int_\Omega J(|x-y|)m(x)n(y){\rm d}x{\rm d}y\ .\label{free2} \end{eqnarray}
Here, $J$ is a non negative, decreasing and  compactly supported function  on $\R_+$ with range $R$.
Notice that we must impose more conditions on $J$ in the case of two species that we did in the single component model.  The reasons for this will be made clear  in Section 6.

The problem considered in  \cite{CCELM} is to minimize ${\cal F}(m,n)$ subject to the  constraint that
\begin{equation}\frac{1}{|\Omega|}\int_\Omega m(x){\rm d}x  \qquad{\rm and}\qquad 
\frac{1}{|\Omega|}\int_\Omega n(x){\rm d}x \label{eq52}
\end{equation}
have certain prescribed values.  As shown in    \cite{CCELM},
this system undergoes a segregating phase transition when $beta$ is large enough for the interaction term to overcome the entropy terms in ${\cal F}$. These would prefer to have $m$ and $n$ to be uniform and this will indeed be the minimizing state for small $\beta$, i.e. high temperature $\beta^{-1}$.
However, for large values of $\beta$,
the advantages of segregation can dominate, and the  fluid separates into two phases, one rich in particles of species $1$, and the other rich
in particles of species $2$.  Our concern here is with the profiles of the densities at  the interface between the two phases.

The nature of the two phases in the bulk is determined by considering the zero range model, in which the length scale  $R$ of the interaction $J$ 
is negligible compared to the size of $\Omega$. Formally this corresponds to setting $J(x-y)=\widehat J\delta(x-y)$. It is also convenient to drop the constraint (\ref{eq52}) and  to consider the function
\begin{equation}f_{\beta, \lambda_1,\lambda_2}(m,n)=m\ln m +n\ln n +\beta \hat J mn -\lambda_1 m-\lambda_2 n\label{local2spec}
\end{equation}
as a local free energy density. Here, as in the one component case, 
$\widehat J = \int_{\R^n} J(|x|){\rm d}x\ $, and
$\lambda_1$ and $\lambda_2$ are Lagrange multipliers that ensure the constraint (\ref{eq52}) on the
total particle numbers. One may also think of $\lambda_1$ and $\lambda_2$ as specified chemical potentials and then determine $m$ and $n$ as functions of $\lambda_1$ and $\lambda_2$.

In \cite{CCELM} it is proved that, if $\lambda_1\neq \lambda_2$, there is an unique couple $\bar m,\bar n$ minimizing $f_{\beta,\lambda_1,\lambda_2}(m,n)$. However, if $\lambda_1=\lambda_2$, there is a $\beta_c$ such that, if $\beta\le\beta_c$ the minimizer is still unique, while, if $\beta>\beta_c$ there are densities $\rho^-<\rho^+$ such that the couples $(\rho^+,\rho^-)$ and $(\rho^-,\rho^+)$ are both minimizers of  $f_{\beta, \lambda_1,\lambda_2}(m,n)$. We focus on the last case.   Thus, in what follows $\lambda_1=\lambda_2=\lambda$. Then $f_{\beta, \lambda_1,\lambda_2}$ is the local Gibbs free energy density.

Analysis of the zero range model suffices to determine the quantity of the fluid that is present in each phase, but not
the surface tension across the boundary.   
We now turn to the  variational problem that determines the density profiles across the interface, and the surface tension. We will assume that the geometry of $\Omega$ is such that the interface is perpendicular to the first coordinate axis; e.g., we take $\Omega$ to be a very long cylinder along the $x_1$-axis with periodic boundary conditions along the other coordinate axes.

First, we need to introduce the one dimensional version of $J$. Choose coordinate $(s,t)$ on $\R^n$ with $s\in \R$ and $t\in \R^{n-1}$, and define $\bar J$ on $\R$ by
$$\bar J(s) = \int_{\R^{n-1}}J(\sqrt{s^2+|t|^2}){\rm d}t\ ,$$
and then 
$\widehat J\ = \int_{\R}\bar J(s){\rm d}s$.
Let  
$g_{\beta, \lambda}=\inf_{m,n\ge 0}f_{\beta, \lambda,\lambda}(m,n)$.
By what has been noted above,
$$g_{\beta, \lambda} =  f_{\beta, \lambda,\lambda}(\rho^-,\rho^+) =   f_{\beta, \lambda,\lambda}(\rho^+,\rho^-)\ .$$

The functional ${\mathcal G}$ defined by
\begin{equation}{\mathcal G}(m,n)=\int_\R \left[ m(x)\ln m(x) + n(x)\ln n(x)+\beta \int_\R  \bar  J(x-y)m(x)n(y){\rm d}y- g_{\beta,\lambda}\right] {\rm d}x \label{eqG}
\end{equation}
is the {\it excess free energy} of a front. We look for the minimizers of this functional for $\beta>\beta_c$.
The minimum value gives the surface tension across the planar phase boundary. Note that we have let $\Omega\to \R$ and that ${\mathcal G}$ is the free energy per unit $(d-1)$- dimensional area.

Our goal in the next sections is to prove a strict displacement convexity property of  this excess free energy functional, and to show, as a consequence, the uniqueness of the minimizing fronts up to translation.   As in the one component case, a rearrangement inequality will enable us to restrict our attention to monotone profiles.  Let
${\mathcal M}_{\rho^-,\rho^+}$ be the subset of ${\mathcal C}_{\rho^-,\rho^+}$  consisting of monotone increasing profiles,
let  ${\mathcal M}_{\rho^+,\rho^-}$ be the subset of ${\mathcal C}_{\rho^+,\rho^-}$  consisting of monotone decreasing profiles

Our main goal mathematically in what follows is to show that
the functional
$$(m,n) \mapsto 
\int_\R\left[\int_\R  \bar J(x-y)m(x)n(y){\rm d}y-\hat J\rho^+\rho^-\right] {\rm d}x$$
is displacement convex on   ${\mathcal M}_{\rho^-,\rho^+}\times {\mathcal M}_{\rho^+,\rho^-}$, where now we have both an increasing and a decreasing density profile.

We shall prove the displacement convexity results in the next section. This time, we shall require certain moment conditions
to obtain the displacement convexity. Hence, before we can apply these results, we need to show {\it a priori}
that all minimizers have good localization properties. We do this by an analysis of the Euler Lagrange equation.

\subsection{Convexity of the interaction energy for ${\cal G}$} 
\medskip

Define the functional
${\cal I}$ on ${\cal M}_{a,b}\times {\cal M}_{c,d}$ by
\begin{equation}\label{free3} 
{\cal I}(m,n) = 
 \int_\R{\rm d}x \left[\int_R \bar J(x-y)m_1(x)m_2(y){\rm d}y- \hat J \hat m(x)\hat n(x)\right]\ .
 \end{equation}
We assume $J$ to be  non negative, even and compactly supported on $\R$ with range $R$. We define 
$\widehat J$ to be the total mass of $J$,
and we define
$$ \widehat m(x)=
\begin{cases} 
 b\, &\text{for $x\ge 0$}\  , \\
 a\, &\text{for $x< 0$}\  
\end{cases} \qquad{\rm and}\qquad 
\widehat n(x)=
\begin{cases} 
 d\, &\text{for $x\ge 0$}\  , \\
 c\, &\text{for $x< 0$}\  .
\end{cases}
$$

Note that in the special case $a = d = \rho^-$ and $b=c = \rho^+$,
$${\cal I}(m,n) = \int_\R\left[\int_\R  \bar J(x-y)m(x)n(y){\rm d}y-\hat J\rho^+\rho^-\right] {\rm d}x\ .$$
Although this special case is all that is needed for our applications here,  we treat the general case
because the small extra  effort yields a broad new class of jointly displacement convex functionals.

The first step in our analysis  is to rewrite ${\cal I}$ as a functional on probability densities. Let the probability densities $\rho_1$ and
$\rho_2$ be defined by
\begin{equation}\label{rows}
m(x)=a+(b-a)\int_{-\infty}^x\rho_1(t)dt; \qquad n(x)=c+(d-c)\int_{-\infty}^x\rho_2(t)dt\ .
\end{equation}
We the  rewrite the functional in terms of $\rho_1$ and $\rho_2$, and integrate by parts. Formally, one moves an antiderivative from each of $\rho_1$ and $\rho_2$
over to $\bar J$. Since $\bar J$ is positive, integrating it twice produces a convex function $W$, different from the one constructed in the one-component case. This is indeed what happens, but
one must be careful about the boundary terms.  The boundary terms do not vanish, but as we shall see, they depend on the densities
in a very nice way, and altogether, one obtains the desired displacement convexity.

To carry  out this analysis, 
define 
\begin{equation}\label{wdef}
W(x) =
\begin{cases} 
\displaystyle{\int_0^x\left(\int_0^t \bar J(s){\rm d}s\right){\rm d}t}\, 
&\text{for $x>0$}\  , \\
W(-x) &\text{for $x<0$} \ . 
\end{cases}
\end{equation}
Then $W''(x) = \bar J(x)$, $W(0) = 0$, and $W$ is an even convex function. Furthermore,
\begin{equation}\label{1b}
\lim_{x\to\infty}W'(x) =\frac{\hat J}{2}, \quad W(x)=\alpha+ \frac{\hat J}{2}|x| \text{ for } |x|\ge R\ .
\end{equation}

\medskip
\begin{lm} \label{thm1bcn}
Let $m\in {\cal M}_{a,b}$  and 
$n\in {\cal M}_{c,d}$. Let $\rho_1$ and $\rho_2$ be the corresponding probability densities  defined in (\ref{rows}).
Then, provided $\rho_1$ and $\rho_2$ have finite first moments, and with $W$ and $\alpha$ defined as above, 
\begin{eqnarray}
\mathcal{I}(m_1,m_2) &=& 
 (a-b)(d-c)\int_\R\int_\R   W(x-y)\rho_1(x)\rho_2(y)  {\rm d}x{\rm d}y\nonumber\\
 &+&[2(b-a)(d-c)+bc+ad]\alpha\nonumber\\
&- &\frac{\hat J}{2}\int_\R\int_\R x\     \big[(b+a)(d-c)\rho_2(x) + (b-a)(c+d)          \rho_1(x)\big ]{\rm d}x\ .\nonumber\\ \label{FunI}
\end{eqnarray}
 \end{lm}
 \medskip
Note that  $(a-b)(d-c)>0$ for $b>a$ and  $c>d$, which is the case when $a=d = \rho^-$ and $c= b = \rho^+$. Thus, $(a-b)(d-c)W(z)$ is a convex function of $z$ on all of $\R$. It follows in the usual way that the first term on the right is displacement convex.  Since $W$ is strictly convex on the support of $J$, it follows as in the proof of Theorem 
\ref{interconv} that this part of the functional (\ref{FunI}) is in fact strictly convex apart from translation. 
The  second term on the right of (\ref{FunI}) is a constant. The third term is a linear combination
of the first moments of $\rho_1$ and $\rho_2$. Since these first moments are displacement affine, we see that
altogether, $\mathcal{I}(m,n)$ is strictly displacement convex, apart from translation. 

The fact that  Lemma \ref{thm1bcn}  requires a conditions on first moments, while Theorem \ref{interconv} does not, means that
it will be a little more work  to apply Lemma \ref{thm1bcn}:  We shall need an {\it a priori} estimate guaranteeing that for any critical point $(m,n)$
of ${\cal G}$, the corresponding densities have finite first moments.  We shall return to this after first proving the theorem.
\medskip
 
 \noindent{\bf Proof:} 
 We start by considering 
 the integral in $x$ first, on a bounded interval $[-L,L]$. Since  $\displaystyle{\bar J(x-y)=-\frac{\partial^2}{ \partial x\partial y}W(x-y)}$ we have that
\begin{eqnarray}
&\phantom{*}&-\int_{-L}^{L} \frac{\partial^2}{ \partial x\partial y}W(x-y)m(x) {\rm d}x = 
\int_{-L}^{L}\frac{\partial }{ \partial y}W(x-y)(b-a)\rho_1(x)  {\rm d}x  \nonumber\\ 
&\phantom{*}&-\frac{\partial }{ \partial y}W(L-y)m(L)  +
\frac{\partial }{ \partial y}W(-L-y)m(-L) \nonumber\\ 
\end{eqnarray} 
Moreover,
\begin{eqnarray}
&\phantom{*}&\int_{-L}^{L}\int_{-L}^{L}\bar J(x-y)m(x)n(y){\rm d}y =\int_{-L}^{L}\int_{-L}^{L}\frac{\partial }{ \partial y}W(x-y)(b-a)\rho_1(x)n(y)  {\rm d}x{\rm d}y\ +
\nonumber\\ 
&\phantom{*}& \int_{-L}^{L}\left[-\frac{\partial }{ \partial y}W(L-y)m(L)  +
\frac{\partial }{ \partial y}W(-L-y)m(-L) \right] n(y){\rm d}y
\nonumber\\
\end{eqnarray}

Now we integrate by parts once more, this time in $y$:

\begin{eqnarray}
&\phantom{*}&\int_{-L}^{L} \frac{\partial}{ \partial y}W(x-y)n(y)  {\rm d}y= -\int_{-L}^{L} W(x-y)(d-c)\rho_2(y)  {\rm d}y\nonumber\\
&\phantom{*}&\displaystyle{+W(x-L)n(L)  - W(x+L)n(-L)}\ .  \nonumber\\
\end{eqnarray} 
Summarizing,
\begin{eqnarray}
&\phantom{*}&\int_{-L}^{L}\int_{-L}^{L}\bar J(x-y)m(x)n(y){\rm d}y\nonumber \\&&=-(b-a)(d-c)\int_{-L}^{L}\int_{-L}^{L}W(x-y)\rho_1(x) \rho_2(y) {\rm d}x{\rm d}y\ 
\nonumber\\ 
&\phantom{*}& +\int_{-L}^{L}\Big[-\frac{\partial }{ \partial y}W(L-y)m(L)+\frac{\partial }{ \partial y}W(-L-y)m(-L) \Big] n(y){\rm d}y \nonumber\\ 
 &\phantom{*}& +(b-a)\int_{-L}^{L}\left[     W(x-L)n(L)-W(x+L)n(-L)            \right]\rho_1(x) {\rm d}x
\nonumber\\
\end{eqnarray} 
Let us examine the boundary terms
$$B_1:=\int_{-L}^{L}\left[-\frac{\partial }{ \partial y}W(L-y)m(L)+\frac{\partial }{ \partial y}W(-L-y)m(-L) \right] n(y){\rm d}y\ ,$$
$$ B_2:=(b-a)\int_{-L}^{L}\left[     W(x-L)n(L)-W(x+L)n(-L)            \right]\rho_1(x) {\rm d}x$$
We have
$$B_1=m(L)\int_{-L}^{L}W(L-y)(d-c)\rho_2(y) {\rm d}y -m(-L)\int_{-L}^{L}W(-L-y)(d-c)\rho_2(y){\rm d}y$$
$$+m(L)\left[-W(L-y)n(y)\right] _{-L}^{+L} +m(-L)\left[W(-L-y)n(y)\right] _{-L}^{+L} \ .$$
$$=(d-c)\int_{-L}^{L}\left[m(L)W(L-y) -m(-L)W(-L-y)\right]\rho_2(y){\rm d}y$$
$$+m(L)\left[-W(0)n(L)+W(2L)n(-L)\right]+m(-L)\left[W(2L)n(L)-W(0)n(-L)\right] $$
For $2L>R$, where $R$ is the range of the interaction $\bar J$,  the last two terms give
$$(bc+ad)({\hat J}L+\alpha)+{\mathcal O}(1)$$
To compute the other term, we consider, for a function $f$ rapidly decaying,
$\int_{-L}^{L}f(x)W(x+L){\rm d}x$ and $\int_{-L}^{L}f(x)W(x-L){\rm d}x$. We have
$$\int_{-L}^{L}f(x)W(x+L){\rm d}x=\int_{0}^{2L}f(z-L)W(z){\rm d}z=\int_{0}^{R}f(z-L)W(z){\rm d}z +\int_{R}^{2L}f(z-L)(\frac{\hat J}{2}z+\alpha){\rm d}z$$
The first term vanishes in the limit $L\to\infty$ because of the decay of $f$ and of the boundedness of $W(z)$ for $z\in[0,R]$.
The second term becomes, if $\int_{\R}|x|f(x){\rm d}x<\infty$,
$$\int_{R-L}^{L}f(x)(\frac{\hat J}{2}(x+L)+\alpha){\rm d}x=\frac{\hat J}{2}\int_{\R}xf(x){\rm d}x+(\alpha+\frac{\hat J}{2}L)\int_{\R}f(x)\d x+{\mathcal O}(1)$$
In conclusion,
$$\int_{-L}^{L}f(x)W(x\pm L){\rm d}x=\pm \frac{\hat J}{2}\int_{\R}xf(x){\rm d}x+(\alpha +L\frac{\hat J}{2})\int_{\R}f(x)\d x+{\mathcal O}(1)$$

Now we apply this result to $B_2$, where the decaying function is $\rho_1$, to get 
$$B_2=(b-a)\left[-(c+d)\frac{\hat J}{2}\int_{\R}x\rho_1(x){\rm d}x+\alpha(d-c)\int_{\R}\rho_1(x){\rm d}x +\frac{\hat J}{2}L(d-c)\int_{\R}\rho_1(x){\rm d}x\right]+{\mathcal O}(1)$$

Now we apply to $B_1$:

\begin{eqnarray}
& &B_1=(d-c)\left[-(b+a)\frac{\hat J}{2}\int_{\R}x\rho_2(x){\rm d}x+\right. \\
& &\left.\alpha(b-a)\int_{\R}\rho_2(x){\rm d}x +\frac{\hat J}{2}L(b-a)\int_{\R}\rho_2(x){\rm d}x\right] +(bc+ad)({\hat J}L+\alpha)+{\mathcal O}(1)\nonumber 
\end{eqnarray}

Finally,
$$B_1+B_2- \hat J \int_{\R}\hat m(x)\hat n(x) {\rm d}x=[2(b-a)(d-c)+bc+ad]\alpha $$ $$-\frac{\hat J}{2}(b+a)(d-c)\int_{\R}y\rho_2(y){\rm d}y -  \frac{\hat J}{2} (b-a)(c+d)          \int_{\R} x\rho_1(x){\rm d}x+{\mathcal O}(1)$$
\qed
 
 \medskip
 
Lemma \ref{thm1bcn} is the key ingredient to prove the analog of Theorem \ref{interconv} for the two-component model introduced in the beginning of this section. We 
now return to this model, and shall apply the lemma with 
 $a=d=\rho^-$ and $b=c=\rho^+$. Let $(w_1,w_2)$ and $(v_1,v_2)$ be in ${\mathcal M}_{\rho^-,\rho^+}\times {\mathcal M}_{\rho^+,\rho^-}$, with corresponding probability densities $(\eta_1,\eta_2)$ and $(\zeta_1,\zeta_2)$, and let $T_1$, $T_2$ be  the monotone maps such that $\zeta_i=T_i\#\eta_i$, $i=1,2$. Moreover, let $\lambda \mapsto (m_\lambda,n_\lambda)$ be the displacement interpolations between $(w_1,w_2)$ and $(v_1,v_2)$ and $T_i^\lambda(x)=\lambda x+(1-\lambda)T_i(x)$.
 
 \medskip
\begin{theorem}\label{interconvmult} 
Suppose that the probability densities $\eta_i$ and $\zeta_i$, $i=1,2$ have finite first moments. Then for $0 < \lambda < 1$,
 \begin{eqnarray}\label{conv1mult}
 &&\mathcal{G}(m_\lambda,n_\lambda)\le (1-\lambda)\mathcal{G}(w_1,w_2)+\lambda \mathcal{G}(v_1,v_2)\ .\nonumber\\
\end{eqnarray}
If  $J$ is strictly positive on some interval, and $(w_1,w_2)$ have  strictly positive derivatives almost everywhere, there is equality if and only if $(v_1,v_2)$ is a translate of $(w_1,w_2)$.
\end{theorem}
\medskip
\noindent{\bf Proof.}
Lemma \ref{thm1bcn} is applicable by the assumption that the probability densities have finite first moments. We set $S_i(x)=T_i(x)-x$, so that 
$${\cal I}(m_\lambda,n_\lambda) = (\rho^+-\rho^-)^2\int_\R\int_\R   W[(x-y)+\lambda(S_1(x)-S_2(y))]d\eta(x)d\eta(y)+\mathcal{A}(m_\lambda,n_\lambda)\ ,$$
with $\mathcal{A}$ affine. The function $W$ is convex on all $\mathbb{R}$, thus the interaction part of the ${\cal G}$ is strictly displacement convex. 
Then remaining term  is simply a linear combination of functions of $m$ and $n$ to which we can apply Theorem \ref{affine}. The strict convexity up to translations follows as in the proof of Theorem \ref{interconv}. \qed

\medskip

\noindent{\bf Remark:} In the two component case we need to use two monotone maps instead of one as in Theorem \ref{interconv}. Therefore it is crucial that $W$ is convex on all of $\mathbb{R}$ and not just
on $(0,+\infty)$ and $(-\infty, 0)$ as in the one component case.
\medskip
 
 We close this section with a corollary showing that one could also use Lemma \ref{thm1bcn} to prove displacement convexity of the interaction energy in the one component model. In fact, in this application,
 the first moment condition drops out.

\begin{cl}\label{altap}  Let $\bar J$ satisfy the conditions below (\ref{free3}), and 
and $W$ defined as in the  (\ref{wdef}).
Let $m$  be  a function that increases monotonically from $-m_\beta$ to $m_\beta$.  Let $\rho$  denote  $m'$, the derivative of $m$. Consider the functional $\Phi(m)$ given by
$\Phi(m) = \int_\R\int_\R \bar J(x-y)\Big[m(x)m(y)  -m^2_\beta\Big]{\rm d}x{\rm d}y$. Then
$$\Phi(m) = 
 -4m_\beta^2\int_\R\int_\R   W(x-y)\rho(x)\rho(y)  {\rm d}x{\rm d}y -6\alpha m_\beta^2\ .$$

\end{cl}

\noindent {\bf Proof}. The functional $\Phi(m)$ is equal to $-\mathcal{G}(m_1,m_2)$ by putting $m_1(x)=m(x)$ and $m_2(x)=-m(x)$.
This shows that $-\Phi$
 is strictly displacement convex, up to translation. \qed

 \medskip

\section{Properties of the minimizers of ${\cal G}$.}
\
\medskip
We restrict our attention to the case  $a=d$, $b=c$.  We need two results on the minimizers for
${\cal G}$, the first of which allows us to restrict our attention to monotone profiles when seeking to minimize ${\cal G}$.
The second guarantees the existence of moments for the two densities corresponding to any minimizing pair $(m,n)$.
These theorems are:

\medskip

\begin{theorem}\label{rear} Suppose that $J(x)$ is even non negative and decreasing. Then any minimizer $(m_1,m_2)$ of ${\cal G}(m_1,m_2)$ in ${\mathcal C}_{\rho^-,\rho^+}\times {\mathcal C}_{\rho^+,\rho^-}$
is monotone in the sense that $m_1$ is increasing and $m_2$ is decreasing.
\end{theorem}

\medskip

This theorem makes it easy to establish the existence of minimizers for {\cal G}. The minimizers satisfy an Euler--Lagrange equation from which we can deduce {\it a priori} estimated needed to apply Lemma \ref{thm1bcn}.

\begin{theorem}\label{inf} Suppose that $J(x)$ is even non negative and decreasing on $\R_+$.
Any  minimizer $w=(w_1,w_2)$ of ${\mathcal G}$ in ${\mathcal C}_{\rho^-,\rho^+}\times {\mathcal C}_{\rho^+,\rho^-}$ 
satisfies
$\rho^-<w_i(x)<\rho^+$
for any $x\in \R$. It has derivative almost everywhere which is strictly positive and with  $\|w_i'\|_{L^1(\R)}$ is bounded. Furthermore, it satisfies the Euler-Lagrange equations
\begin{equation}\label{e-l}\ln m(x)+\beta (J*n)(x)=\mu,\quad \ln n(x)+\beta (J*m)(x)=\mu,\end{equation}
where $\mu=\mu_1-1$ and $*$ denotes convolution.  Its derivative $w$ satisfies almost everywhere the equations
\begin{equation}\label{e-l'}\frac{w_1'(x)}{w_1(x)}+\beta(J*w_2')(x)=0, \quad \frac{w_2'(x)}{w_2(x)}+\beta(J*w_1')(x)=0
\end{equation}
Finally, it converges to its asymptotic values exponentially fast, in the sense that there is $\alpha>0$ such that $(w_1(x)-\rho_\mp)e^{\alpha|x|}\to 0 \text{ as } x\to \mp \infty$ and $
(w_2(x)-\rho_\pm)e^{\alpha|x|}\to 0 \text{ as } x\to \mp \infty$.
\end{theorem}

\medskip

The proof of Theorem \ref{rear} is adapted from a related result in \cite{CCELM} for functions on the $d$ dimensional torus.  One could instead adapt the proof of Alberti's rearrangement inequality
in \cite{alb} and remove the requirement that $J$ be decreasing. But the present approach has the advantage of working also on the torus, and not only the line. 
The proof of the final part of Theorem \ref{inf}, which is important for our application here since it provides the existence of moments, 
 is adapted from the proof of a similar result for the one component system in \cite{dopt2}.   In the rest of this section, we present these proofs.
 
 \medskip

\noindent{\bf Proof of Theorem \ref{rear}:}
To show this, we use a rearrangement inequality similar to those introduced in \cite{CCELM} for the analogous problem in the $d$-dimensional torus. For any $x_0\in\mathbb{R}$, let $T_{x_0}$ denote the reflection about $x_0$:
$$T_{x_0}(x) = 2x_0 -x\ .$$
Then define  $\cal D$, as the set of functions on $\mathbb{R}$ having finite limits at $\pm \infty$  and 
the operators $R^\pm_{x_0}$ on ${\cal D}$ by
\begin{equation}R^+_{x_0}g(x) =
\begin{cases}
      \max\{g(x),g(T_{x_0})\} & \text{ if $x\ge x_0$ }, \\
      \min\{g(x),g(T_{x_0})\} & \text{ if $x\le x_0$}.
\end{cases}
\end{equation}
\begin{equation}R^-_{x_0}h(x) =
\begin{cases}
      \max\{h(x),h(T_{x_0})\} & \text{ if $x\le x_0$ }, \\
      \min\{h(x),h(T_{x_0})\} & \text{ if $x\ge x_0$}.
\end{cases}
\end{equation}
Let us also define
${\displaystyle
\hat g(x)=\begin{cases} 
\displaystyle{\lim_{x\to -\infty} g(x) }&\text{if $x<0$}\\
\displaystyle{\lim_{x\to +\infty}g(x)} &\text{ if $x\ge 0$}
\end{cases}
}$
and $\hat h$ similarly.

For any fixed $x_0$ and $g,h \in \cal D$, let $g^{\star}$ denote 
$R^+_{x_0}g$ and  $h_{\star}= R^-_{x_0}h.$
We now wish to show that 
$$\int_{\mathbb{R}}\left[ \int_{\mathbb{R}} g(x)J(x-y)h(y){\rm d}y -\hat J\hat g(x)\hat h(x)\right]{\rm d}x\ge
\int_{\mathbb{R}}\left[ \int_{\mathbb{R}} g^\star (x)J(x-y) h_\star (y){\rm d}y -\hat J\hat g(x)\hat h(x)\right]{\rm d}x$$
with equality if and only if $g = g^\star$ and $h = h^\star$.

To do this, let $\mathbb{H}_+$ denote the half line $\{x\ |\ x > x_0\}$, and $\mathbb{H}_-$ denote the half line $\{x\ |\ x < x_0\}$ and
observe that
\begin{eqnarray}
\int_{\mathbb{R}} \int_{\mathbb{R}} g(x)J(x-y)h(y){\rm d}x{\rm d}y &  = &\nonumber\\
\int_{\mathbb{H}_+} \int_{\mathbb{H}_+} g(x)J(x-y)h(y){\rm d}x{\rm d}y +
\int_{\mathbb{H}_-} \int_{\mathbb{H}_-} g(x)J(x-y)h(y){\rm d}x{\rm d}y & + &\nonumber\\
\int_{\mathbb{H}_-} \int_{\mathbb{H}_+} g(x)J(x-y)h(y){\rm d}x{\rm d}y +
\int_{\mathbb{H}_+} \int_{\mathbb{H}_-} g(x)J(x-y)h(y){\rm d}x{\rm d}y & = &\nonumber\\
\int_{\mathbb{H}_+} \int_{\mathbb{H}_+} g(x)J(x-y)h(y){\rm d}x{\rm d}y +
\int_{\mathbb{H}_+} \int_{\mathbb{H}_+} g(T_{x_0}x)J(x-y)h(T_{x_0}y){\rm d}x{\rm d}y & + & \nonumber\\
\int_{\mathbb{H}_+} \int_{\mathbb{H}_+} g(T_{x_0}x)J(T_{x_0}x-y)h(y){\rm d}x{\rm d}y +
\int_{\mathbb{H}_+} \int_{\mathbb{H}_+} g(x)J(x-T_{x_0}y)h(T_{x_0}y){\rm d}x{\rm d}y\nonumber\\
\end{eqnarray}

The desired inequality is then a consequence of the following inequality
for pairs of real
numbers: Let $a_1$ and $a_2$ and $b_1$ and $b_2$ be any four positive
real numbers.
Rearrange $a_1$ and $a_2$ to decrease, and $b_1$ and $b_2$ to increase;
i.e., let 
$a^\star_1 = \max\{a_1,a_2\}$, $a^\star_2 = \min\{a_1,a_2\}$, 
$b^\star_1 = \min\{b_1,b_2\}$ and  $b^\star_2 = \max\{b_1,b_2\}$. Then
\begin{equation}
\label{hardy}
a_1^\star b_1^\star + a_2^\star b_2^\star - a_1b_1 - a_2 b_2=
\Delta\ \le 0,
\end{equation}
\begin{equation}
\label{hardy1}
a_1^\star b_2^\star + a_2^\star b_1^\star - a_1b_2 - a_2 b_1=-
\Delta\ \ge 0,
\end{equation}
and there is equality if and only if $a_1 = a_1^\star$ and $b_1 =
b_1^\star$ or
$a_1^\star=a_2$ and $b_1^\star=b_2$.

We now apply the above inequalities with
\begin{equation}
\label{A.15}
a_1 = g(x )\quad a_2 = g(T_{x_0}x)\quad b_1 =
h(y)\quad{\rm and} \quad
b_2 = h(T_{x_0}y)\ .
\end{equation}
Then
\begin{equation}
\label{A.16}
a_1^\star = R_{x_0}^+g(x)\quad a_2^\star = R_{x_0}^+g(T_{x_0}x)
\quad b_1^\star = R_{x_0}^-h(y)\quad{\rm and}\quad b_2^\star =
R_{x_0}^-h(T_{x_0}y)\ .
\end{equation}
Since
$J(T_{x_0}x-y)=J(x-T_{x_0}y)< J(x-y)$, we get
\begin{eqnarray} \label{aaa}
g(x)J(x -y)h(y) +g(T_{x_0}x)J(x -y)h(T_{x_0}y) &+&\nonumber\\
g(T_{x_0}x)J(T_{x_0}x -y)h(y)  +   g(x)J(T_{x_0}x -y)h(T_{x_0}y)  &-&\nonumber\\
R_{x_0}^+g(x)J(x -y)R_{x_0}^-h(T_{x_0}y) +R_{x_0}^+g(T_{x_0}x)J(x -y)R_{x_0}^-h(T_{x_0}y) &-&\nonumber\\
R_{x_0}^+g(T_{x_0}x)J(T_{x_0}x -y)R_{x_0}^-h(y)  + 
R_{x_0}^+g(x)J(T_{x_0}x -y)R_{x_0}^-h(T_{x_0}y) &=&\nonumber\\
-
\Delta\big[J(x -y)-J(x-T_{x_0}y)\big]\ge &0&\nonumber\\
\end{eqnarray}
for almost every $x$ and $y$ in $\mathbb{H}_+$, with equality if
and only if
\begin{equation}
\label{AA11}
g(T_{x_0}x) \le g(x)\qquad{\rm and}\qquad h(T_{x_0}y)
\ge h(y)
\end{equation}
or
\begin{equation}
\label{AA22}
g(T_{x_0}x) \ge g(x)\qquad{\rm and}\qquad h(T_{x_0}y)
\le h(y)
\end{equation}
for almost every $x$ and $y$ in $\mathbb{H}_+$. Now unless $g$ is
constant, we can find 
$x$ and  $x_0$ so that either $g(T_{x_0}x) < g(x)$ or
$g(T_{x_0}x) > g(x)$.
Suppose it is the first case. Then (\ref{AA11}) holds, and
for almost every $y$, we must have $h(T_{x_0}y)
\ge g(y)$. Making a similar argument for $h$, we see that one of
(\ref{AA11}) or (\ref{AA22})
must hold for almost every $x$ and $y$.
The only way that
this can happen is if $g$ and $h$
are monotone. Now, by integrating (\ref{aaa}) on $\mathbb{H}_+$ we conclude the proof. \qed

\medskip

 \medskip

\noindent{\bf Proof of Theorem \ref{inf}:}  Everything but the exponential decay is standard, and details of the proofs of
similar results can be found in \cite{CCELM}. To prove the exponential decay, we once again take advantage of the finite range $R$ of $J$. 

Define a transformation $\Phi: \R^2 \to \R^2$ by
$\Phi(m,n) = (e^{\mu - \beta \widehat J n} ,  e^{\mu - \beta \widehat J m})$.
Then $(\rho^+,\rho^-)$ and $(\rho^-,\rho^+)$ are two stable fixed points of $\Phi$; the Jacobian of $\Phi$, $D\Phi$,
is a strict contraction at either of them.  
Thus, there is a $\delta>0$ and an $\epsilon>0$ so that if 
$$|m-\rho^+| + |n-\rho^-| < \delta     \quad\Rightarrow \quad
\|D\phi(m,n)\| < 1- \epsilon\ .$$ 

Now, consider any minimizer $w = (w_1,w_2)$ with $\lim_{x\to\infty}w_1(x) = \rho^+$ and 
$\lim_{x\to\infty}w_2(-x) = \rho^-$.  Then there is an $L<\infty$ so that
$x \ge L \quad\Rightarrow \quad  |w_1(x)-\rho^+| + |w_2(-x)-\rho^-| < \delta$.
Now for $x> L+R$, 
$$\frac{J}{\widehat J}*w_1(x) \ge \rho^+-\delta \qquad{\rm and}\qquad  
\frac{J}{\widehat J}*w_2(x) \le \rho^- +\delta\ .$$
Since
${\displaystyle (w_1(x),w_2(x)) = \Phi\left(\frac{J}{\widehat J}*w_1(x), \frac{J}{\widehat J}*w_2(x)\right)}$,
it follows that
for   $x> L+R$, 
$ |w_1(x)-\rho^+| + |w_2(-x)-\rho^-| < (1-\epsilon)\delta$.
Iterating this argument leads to the conclusion that  for   $x> L+kR$, 
$ |w_1(x)-\rho^+| + |w_2(-x)-\rho^-| < (1-\epsilon)^k\delta$.
A similar argument applies as $x$ tends to $-\infty$. \qed
\medskip

\section{For the functional ${\cal G}$, critical points are minimizers}

We are now ready to prove the main theorem for $\G$:

\medskip

\begin{theorem} If $(w_1,w_2)$ and  $(v_1,v_2) $ are any two critical  points of ${\cal G}$ in ${\mathcal M}_{\rho^-,\rho^+}\times {\mathcal M}_{\rho^+,\rho^-}$, then there is an $a\in R$ so that
\begin{equation}\label{trans}
(v_1(x),v_2(x))  =  (w_1(x-a),w_2(x-a))\ .
\end{equation}
Thus, there is exactly one critical point  $(w_1,w_2)$  such that $w_1(0) =w_2(0)$. It is symmetric in the sense that
$w_1(x) = w_2(-x)$ for all $x$. 
\end{theorem}

\noindent{\bf Proof:} 
We keep the notation of Section 5. Theorem \ref{interconvmult} is applicable since by Theorem \ref{inf}, the probability densities $\eta_i$ and $\zeta_i$, $i=1,m$, have finite moments of every order. Now, if $(m_\lambda,n_\lambda)$ is the displacement convex interpolation between
$(w_1,w_2)$ and  $(v_1,v_2) $, ${\cal G}(m_\lambda,n_\lambda)$ is constant since both endpoints are critical points.
By the strict convexity up to translation, we see that (\ref{trans}) is true.

Since $(w_1,w_2) \in  {\mathcal M}_{\rho^-,\rho^+}\times {\mathcal M}_{\rho^+,\rho^-}$, and both functions are strictly monotonic, there is
some $b$ such that $w_1(b) = w_2(b)$.   Because of the strict monotonicity of 
$w_1$; i.e., the strict positivity of its derivative, which was proved in Theorem \ref{inf}, this value of $b$ is unique.

Next, by the symmetries of the functional, since 
$(w_1(x),w_2(x))$  is any minimizer of ${\cal G}$ in ${\mathcal M}_{\rho^-,\rho^+}\times {\mathcal M}_{\rho^+,\rho^-}$,
then so is  $(w_2(-x),w_1(-x))$. Hence, by the first part of the Theorem, there is an $a\in R$ so that
\begin{equation}\label{trans2}
(w_2(-x),w_1(-x)) =  (w_1(x-a),w_2(x-a))\ .
\end{equation}
Evaluating both sides at $x=0$, we see that since  $w_1(0) =w_2(0)$,  $w_1(-a) =w_2(-a)$. By the uniqueness of the crossing point
established above,  $a=0$, so that 
$(w_2(-x),w_1(-x)) =  (w_1(x),w_2(x))$
for all $x$. \qed

\section{Stationary monotone profiles in several dimensions.}

We close the paper by pointing out that our analysis of the two component case can be adapted to yield
a uniqueness theorem for the one component case in higher dimensions. 

Let $\Omega$ be a $(d-1)$-dimensional cube of size $L$ spanned by the orthogonal vectors $e_1, \dots,e_{d-1}$ and ${\cal C}_{a,b, \Omega}$ be the set of continuous functions $m(x,y)$ from $\R\times\R^{d-1} $ to $\R$ such that for all $y\in \R^{d-1}$
$$\lim_{x\to -\infty} m(x,y) =a\qquad{\rm and}\qquad  \lim_{x\to +\infty} m(x,y) =b\ ,$$
and  such that $m$ is $L$-periodic on $R^{d-1}$ in the sense that
$m(x, y+L e_k)=m(x,y)$
for each $k=1,\dots,d-1$ and for each $y\in \R^{d-1}$. 

Consider the following 
 $d$-dimensional free energy on ${\cal C}_{-a,a, \Omega}$
\begin{eqnarray}
\F(m) &=&\int_{\R\times \Omega}F(m(x,y))\d x \d y\\ \nonumber &-& \frac{1}{2}\int_{\R\times \Omega}\int_{\R\times \Omega} (m(x_1,y_1) - m(x_2,y_2))^2J(x_1-x_2, y_1-y_2)\d x_1\d x_2\d y_1\d y_2\ ,
\end{eqnarray}
$J(x,y)=U(\sqrt{x^2+|y|^2})$, with $U$ monotone decreasing, finite range smooth function on $[0,+\infty)$ and $F$ an even double well potential with minima in $-a$ and $a$ and $F(\pm a)=0$. (These specific conditions on $F$ enable us to be brief, and can easily be relaxed.)

Obviously, if $\bar m(x)$ is a minimizer for the corresponding one dimensional problem, then
$$\bar m(x,y) := \bar m(x)$$
is a critical point of $\F$ on ${\cal C}_{-a,a, \Omega}$, and is an obvious candidate to be the unique minimizer.
We shall show here that not only is it the minimizer -- this fact has been proved by Alberti \cite{alb} --
but that, up to translation in $x$, $\bar m(x,y)$ is the {\em only} solution of the Euler-Lagrange equation
for minimization of $\F$ that is monotone in $x$ for all $y$.  
A related  question as to  whether all monotone solutions of the Euler-Lagrange equation have this special form has been extensively investigated for the local variant of the free energy (Allen-Cahn or van der Waals) with $\int|\nabla m(x,y)|^2$ in place of the non-local interaction integral above.  It turns out
that the non-local case may be easily treated by regarding the one dimensional profiles $x\mapsto m(x,y)$
for different $y$ as profiles for different components, and applying our previous results.

Define ${\cal M}_{a,b, \Omega}$ to  be the subset of ${\cal C}_{a,b, \Omega}$ for which $m(x,y)$
is monotone in $x$ for each $y\in \R^{d-1}$. As before, there is a rearrangement inequality that allows one to reduce the minimization problem over $ {\cal C}_{-a,a,\Omega}$ to minimization over ${\cal M}_{-a,a,\Omega}$:
Given $m\in {\cal C}_{-a,a,\Omega}$ we define $m^*\in {\cal M}_{-a,a,\Omega}$ as follows by separately rearranging $m(\,\cdot\, ,y)$ for each $y\in \R^{d-1}$, using the one dimensional rearrangement procedure. By the rearrangement results cited above, $\F(m^*) \le F(m)$. In \cite{alb}, Alberti proceeds with
a careful study of the cases of equality here. Instead, we henceforth restrict our attention to $m \in 
{\cal M}_{-a,a,\Omega}$, and shall show that up to translation in $x$,  there is just one solution of the 
Euler-Lagrange equation in this set.

\begin{eqnarray}
&&\F(m) =\int_{\R\times \Omega}[F(m(x,y))-\hat J( m^2(x,y)-a^2)]\d x\d y+\\ \nonumber && \int_{\R\times \Omega}\left[\int_{\R\times \Omega} m(x_1,y_1)m(x_2,y_2)J(x_1-x_2, y_1-y_2)\d x_1\d y_1- \widehat J a^2\right]\d x_2\d y_2\ ,
\end{eqnarray}
where $\widehat J = \int_{\Omega\times \R}J(x,y)\d x\d y$. 
We now observe that the second term on the right can be written in terms of the ${\cal I}$ functional that has been studied in Section \ref{binsec}. Indeed,
this term can be written as
$$-\int_{\Omega\times\Omega}{\cal I}(m(\,\cdot\, , y_1),
-m(\,\cdot\, , y_2))\d y_1 \d y_2\ .$$
This identity relates the multidimensional problem to the two species problem:
here $m(\,\cdot\, , y_1)$ plays the role of the profile for one species, and
$- m(\,\cdot\, , y_2)$ plays the role of the profile for the other species.

Now notice that for $a+b =0$ (or $c+d =0$), the statement of Lemma \ref{thm1bcn} simplifies in a significant way: The first moments drop out as in Corollary ~\ref{altap}, and we have (using the notation from the lemma)
$$
-\mathcal{I}(m_1,-m_2) = 
4a^2\int_\R\int_\R   W(x_1-x_2)\rho_1(x_1,y_1)\rho_2(x_2,y_2))  {\rm d}x_1{\rm d}x_2 + 6a^2\alpha\ .$$
(First moments could be dealt with as before, but we avoid doing so in order to focus on how one may regard the multidimensional problem as a multi-component problem, which is the main point of this section.)

Given two profiles $m_0$ and $m_1$ in ${\cal M}_{-a,a,\Omega}$, let
$m_\lambda$ be the interpolation defined by interpolating between
$m_0(\,\cdot\, , y)$ and $m_1(\,\cdot\, , y)$ separately in each $y$.
Let $x \mapsto T(x,y)$ be the corresponding optimal transportation plan, and let $S(x,y) = T(x,y) -x$. 
Let $m_\lambda(x,y)$ be the induced interpolation between  $m_0(x,y)$ and $m_1(x,y)$.
Then
\begin{multline}-{\cal I}(m_\lambda(\,\cdot\, , y_1),
-m_\lambda (\,\cdot\, , y_2)) = \\
4a^2\int_\R\int_\R   W[x_1-x_2 + \lambda(S(x_1,y_1) - S(x_2,y_2)]\rho_1(x_1,y_1)\rho_2(x_2,y_2)  {\rm d}x_1{\rm d}x_2 + 6a^2\alpha\ .\nonumber
\end{multline}
Since $W$ is strictly convex near the origin if $J$ is strictly positive near the origin, it follows that
if $y_2$ and $y_1$
sufficiently close to one another,  $\rho_1(x_1,y_1)\d x_1$ and $\rho_2(x_2,y_2)\d x_2$ both assign positive mass to some small interval around some $x_0$. Therefore, for such $y_1$ and $y_2$, we see that $\lambda \mapsto 
-{\cal I}(m_\lambda(\,\cdot\, , y_1),
-m_\lambda (\,\cdot\, , y_2))$ is strictly convex, and for {\em any} $y_1$ and $y_2$ it is convex. Clearly, the set of points $(y_1,y_2)$ for which we have  strict convexity is a set of positive measure
(containing the diagonal)  with respect to $\d y_1\d y_2$, and so 
$$\lambda \mapsto 
-\int_{\Omega\times\Omega}{\cal I}(m_\lambda(\,\cdot\, , y_1),
-m_\lambda(\,\cdot\, , y_2))\d y_1 \d y_2\ $$
is strictly convex, apart from translation in $x$.  This strict convexity proves that, up to translation in $x$,
there is just one critical point of $\F$ in  ${\cal M}_{-a,a,\Omega}$. Since clearly $\bar m (x,y)$ is a critical point, we have the following:

\begin{theorem}\label{manyd}  Assume that $J$ is bounded below by a strictly positive number on some neighborhood of the origin. Let $m(x,y)$ be any solution of the Euler-Lagrange equation for the minimization of $\F$ that belongs to ${\cal M}_{-a,a,\Omega}$.  Then for some $x_0\in \R$,
$m(x,y) = \bar m(x-x_0)$
for al $x$ and $y$, where $\bar m$ the antisymmetric minimizer for one dimension. 
\end{theorem}

\medskip
\noindent{\bf Acknowledgement.} {This work was partially supported by U.S. National Science Foundation grants DMS 06-00037 and DMR 01-279-26, AFOSR grant AF-FA9550-04, MIUR, GNFM-INDAM, and POCI/MAT/61931/2004.  In addition, the authors thank the Erwin Sch\"odinger Institute for hospitality
during June and July 2006 when much of the content of the paper was worked out.}

\end{document}